\documentclass[10pt]{article}
\font\smallit=cmti10

\usepackage{amssymb,amsmath,amsthm} 
\usepackage{graphicx,ctable,booktabs}
\usepackage{amsfonts}
\usepackage[section]{placeins}
\usepackage[mathscr]{euscript}
\usepackage[margin=0.5cm]{caption}
\usepackage{subcaption}
\usepackage{listings}
\newenvironment{definition}[1][Definition]{\begin{trivlist}
\item[\hskip \labelsep {\bfseries #1}]}{\end{trivlist}}

\newcommand{\ZZ}{\mathbb{Z}}
\newcommand{\QQ}{\mathbb{Q}}
\newcommand{\RR}{\mathbb{R}}

\renewcommand{\NN}{\mathbb{N}}

\newtheorem{theorem}{Theorem}
\newtheorem{proposition}{Proposition}
\newtheorem{lemma}{Lemma}

\makeatletter

\renewcommand\section{\@startsection {section}{1}{\z@}
{-30pt \@plus -1ex \@minus -.2ex}
{2.3ex \@plus.2ex}
{\normalfont\normalsize\bfseries}}

\renewcommand\subsection{\@startsection{subsection}{2}{\z@}
{-3.25ex\@plus -1ex \@minus -.2ex}
{1.5ex \@plus .2ex}
{\normalfont\normalsize\bfseries}}

\renewcommand{\@seccntformat}[1]{\csname the#1\endcsname. }

\makeatother

\begin{document}

\begin{center}
\uppercase{Quadratic Packing Polynomials on Sectors of $\RR^2$}
\vskip 20pt
{\bf Madeline Brandt\footnote{Research supported by the National Science Foundation (grant number DMS 1358659) and the National Security Agency (NSA grant H98230-13-1-0273).}}\\
{\smallit Reed College, Portland, OR 97202, United States}\\
{\tt mbrandt@reed.edu}\\
\vskip 30pt
%\centerline{\smallit Received: , Revised: , Accepted: , Published: } % We will fill in the dates
\vskip 30pt
\end{center}

\centerline{\bf Abstract}
\noindent
%A polynomial $p:\RR^2 \rightarrow \RR$ on a set $S \subset \RR^2$ is a packing polynomial if $p|_{\ZZ^2}$ gives a bijection from $S \cap \ZZ^2$ to $\NN$. I provide a necessary and sufficient condition for the existence of quadratic $k$-stair packing polynomials on $S(n/m)$
A polynomial $p(x,y)$ on a region $S$ in the plane is called a packing polynomial if the restriction of $p(x,y)$ to $S\cap \ZZ^2$ yields a bijection to $\NN$. In this paper, we determine all quadratic packing polynomials on rational sectors of $\RR^2$.

%
%\pagestyle{myheadings} 
%\markright{\smalltt INTEGERS: 14 (2014)\hfill} 
%\thispagestyle{empty} 
%\baselineskip=12.875pt 
%\vskip 30pt

%%%% 
%%%%%%%%%%%
\section{Introduction}

Let $S \subseteq \RR^2$, and let $I = S \cap \NN^2$. A polynomial $f : \RR^2 \rightarrow \RR$ is a \emph{packing polynomial on $S$} if $f|_I$ is a bijection from $I$ to $\NN$. In 1923 Fueter and P\'{o}lya \cite{fueterpolya} proved that the \emph{Cantor polynomials},
$$
\begin{array}{ccc}
f(x,y) = \frac{(x+y)^2}{2} + \frac{x + 3 y}{2}& \text{ and } &
g(x,y) = \frac{(x+y)^2}{2} + \frac{3x + y}{2},
\end{array}
$$
are the only quadratic packing polynomials on $\RR_{\geq 0}^2$, and Vsemirnov \cite{vsemirnov} gives two elementary proofs of this theorem.
%The Cantor polynomials enumerate points in $\NN^2$ along lines $y = - x + c$ for integers $c$.
%(see Figure \ref{cantor})
Fueter and P\'{o}lya also conjectured that the Cantor polynomials are in fact the only packing polynomials on $\NN^2$.
 In 1978, Lew and Rosenberg \cite{LR1} showed that there are no cubic or quartic packing polynomials on $\NN^2$, but the existence of higher degree packing polynomials remains unknown.\\

%%%%% CANTOR POLY FIGURE
%\begin{figure}
%\centering
%\begin{subfigure}{.5\textwidth}
%  \centering
%  \includegraphics[width=.8\linewidth]{cantor1}
%  \caption{Values of $f(x,y)$ on $\NN^2$.}
%  \label{fig:sub1}
%\end{subfigure}%
%\begin{subfigure}{.5\textwidth}
%  \centering
%  \includegraphics[width=.8\linewidth]{cantor2}
%  \caption{Values of $g(x,y)$ on $\NN^2$.}
%  \label{fig:sub2}
%\end{subfigure}
%\caption{Behavior of the Cantor polynomials.}
%\label{cantor}
%\label{fig:test}
%\end{figure}

In this paper, we study quadratic packing polynomials on rational sectors. For all $\alpha \in \RR_{\geq 0}$, let
$$
S(\alpha) = \{(x,y) \in \RR^2 \ |\ x,y \geq 0 \text{ and }y\leq \alpha x\},
$$
and let $I(\alpha)$ be the set of lattice points contained in $S(\alpha)$. If $\alpha$ is an integer  (rational, irrational), we call $S(\alpha)$ an integral (rational, irrational) sector. The following results are known for quadratic packing polynomials on rational sectors. \\

In 2013, Nathanson \cite{Nathonson} gave two quadratic packing polynomials on $S(n)$, for $n\in \NN$,
$$
f_n(x,y) = \frac{n}{2}x^2+\left(1-\frac{n}{2}\right)x+y,
$$

$$
g_n(x,y) = \frac{n}{2}x^2+\left(1+\frac{n}{2}\right)x-y.
$$

Subsequently, Stanton \cite{stanton}  proved that these polynomials, along with four polynomials on $S(3)$ and $S(4)$, are the only quadratic packing polynomials on integral sectors.\
After classifying the polynomials on integral sectors, Stanton discovered a necessary condition for quadratic packing polynomials on rational sectors.

\begin{theorem}[Stanton \cite{stanton}]
Let $n/m\geq1$ and $(n,m)=1$. Suppose $S(\frac{n}{m})$ has a quadratic packing polynomial $p$, and let $p_2(x,y)$ denote the homogeneous quadratic part of $p$. Then $n$ divides $(m-1)^2$, and 
$$
p_2(x,y) = \frac{n}{2}\left(x - \frac{m-1}{n}y\right)^2.
$$
\label{csthm}
\end{theorem}

We observe that the restriction $n/m \geq 1$ does not result in any loss of generality because there is a bijection, observed by Nathanson in \cite{Nathonson},
 from $I(n/m)$ to $I\left(\frac{n}{m-n\lfloor m/n \rfloor}\right)$ given by
$$
W_{n/m} = \left (\begin{array}{c c}
1 & -\lfloor m / n \rfloor \\
0 & 1
\end{array} \right ).
$$
In light of this, we will say that two packing polynomials $p$ on $S(\alpha)$ and $q$ on $S(\beta)$ are \emph{equivalent} if there exists a linear map $T:I(\alpha) \rightarrow I(\beta)$ which is a bijection from $I(\alpha)$ to $ I(\beta)$ such that $p = q \circ T$.\\

%%%% RESULTS FROM NATHONSON ABOUT SMALL ALPHA
%Nathanson \cite{Nathonson} gives the following result.
%\begin{theorem} Let $n,m$ be relatively prime positive integers with $1 \leq n < m$, and $n | m-1$. Let $d = (m-1)/n$. Then
%$$
%f_{n/m}(x,y) = \frac{n}{2}(x-d y)^2+\frac{(2-n)x+(dn-2d+2)y}{2},
%$$
%$$
%g_{n/m}(x,y)=\frac{n}{2}(x-d y)^2+\frac{(2+n)x+(2d+m+1)y}{2}
%$$
%are quadratic packing polynomials on $I(r/s)$.
%\label{nathansonthm}
%\end{theorem}
%\begin{remark} 
%Applying Stanton's bijection map also gives all quadratic packing polynomials on sectors $I(n/m)$, where $n<m$ and $n \mid m-1$, because under her bijection $I(n/m)$ becomes an integral sector.
%\end{remark}

In this paper, we determine all quadratic packing polynomials on rational sectors up to equivalence by finding the necessary equations for quadratic packing polynomials on rational sectors, and then by finding a sufficient condition for the resulting polynomials to be packing polynomials.
 In Section \ref{kstairsection}, we start by introducing the notion of a $k$-stair polynomial, giving some basic results on their properties, and demonstrating that all quadratic packing polynomials must be $k$-stair polynomials. 
We proceed to give necessary and sufficient conditions for $k$-stair polynomials to be packing polynomials in Section \ref{nsconditions}.
 We conclude in Section \ref{knot4} with our main result: the classification of all quadratic packing polynomials on rational sectors.

%%%%%%%%%%
%\subsection*{Open Problems}
%Some open problems to consider include:
%
%\begin{enumerate}
%\item Is Stanton's necessary condition sufficient?
%\item Can we characterize all packing polynomials on rational sectors?
%\item Are packing polynomials on irrational sectors?
%\end{enumerate}

%%%%%%%%%%%
\section{$k$-Stair Polynomials}
\label{kstairsection}

For the remainder of this paper, assume that $m$ and $n$ are relatively prime, the integer $n
$ divides $(m-1)^2$, and let $l = (n,m-1)$. Let $p(x,y)$ be a packing polynomial, so that by Theorem \ref{csthm} we may write 
$$p(x,y)= \frac{n}{2}\left(x - \frac{m-1}{n}y\right)^2+dx+ey+f.$$

\begin{definition}
We call the line segment $$y = \frac{n}{m-1}x-c\frac{l}{m-1}$$ for $c \in \NN$ and $(x,y) \in S(n/m)$ the \emph{$c^{th}$ staircase of $I(n/m)$}. A \emph{stair} is a point with integer coordinates on a staircase. The \emph{first} stair on the $c^{th}$ staircase is the stair with minimal $x$-coordinate. Two stairs $r,s$ are \emph{consecutive} if they are on the same staircase and there is no other stair on the line segment from $r$ to $s$. For $c \in \NN$, define
$$
S_c \equiv \left\{(x,y) \in I\left(\frac{n}{m}\right) \ |\ y = \frac{n}{m-1}x-c\frac{l}{m-1}\right\}.
$$
 \end{definition}

\begin{lemma} We have $I(n/m) =\cup_{c \in \NN} S_c$.
\begin{proof} Clearly $I(n/m) \supseteq \cup_{c \in \NN} S_c$. For the other direction, let $(a,b) \in I(\frac{n}{m})$, and let $h =a \frac{n}{l} - b \frac{m-1}{l}$. Consider the following line with slope $\frac{n}{m-1}$ through the point $(a,b)$:
$$
y = \frac{n}{m-1} x - \frac{l}{m-1}h.
$$
Since $l \mid n$ and $l \mid m-1$, and $b/a \leq n/m$, we have $h \in \NN$. Therefore, $(a,b)$ is a stair on $S_{h}$.
\end{proof}
\end{lemma}

\begin{lemma} If $p$ is a quadratic packing polynomial on $S(\frac{n}{m})$, and $(x,y) \in S(\frac{n}{m})$, then for some $k \in \NN$,
$$
p\left(x+\frac{m-1}{l}, y+ \frac{n}{l}\right) - p(x,y) = \pm k.
$$
\begin{proof}
By Stanton's necessary condition, $p_2(x,y) = \frac{n}{2}(x-\frac{m-1}{n}y)^2$. If $L$ is a staircase, then $p|_L$ is linear because $p_2|_L$ is constant.
%If $(x,y)$ is a stair, then $(x+\frac{m-1}{l}, y+ \frac{n}{l})$ and $(x,y)$ are consecutive stairs, and so $p(x,y) \in \NN$ and $p(x+\frac{m-1}{l}, y+ \frac{n}{l}) \in \NN$. Thus, $p(x,y)-p(x+\frac{m-1}{l}, y+ \frac{n}{l}) \in \ZZ$, and since $p|_L$ is linear, this will be a fixed constant $k$ for arbitrary $(x,y)$. 
\end{proof}
\label{restriction1}
\end{lemma}

\begin{definition}
Let $p:S(\frac{n}{m}) \rightarrow \RR$ be a quadratic polynomial with $p_2(x,y) = \frac{n}{2}(x-\frac{n}{m-1}y)^2$ and $p(\ZZ^2)\subseteq \ZZ$. Then $p$ is a \emph{$k$-stair} polynomial if for any two consecutive stairs $r,s$, we have $p(r)-p(s) = \pm k$. If $|r| < |s|$ and $p(s)-p(r)=k$, then $p$ we call \emph{ascending}, otherwise we call $p$ \emph{descending}.
\end{definition}

\begin{figure}
\centering
\begin{subfigure}{.5\textwidth}
  \centering
  \includegraphics[width=.9\linewidth]{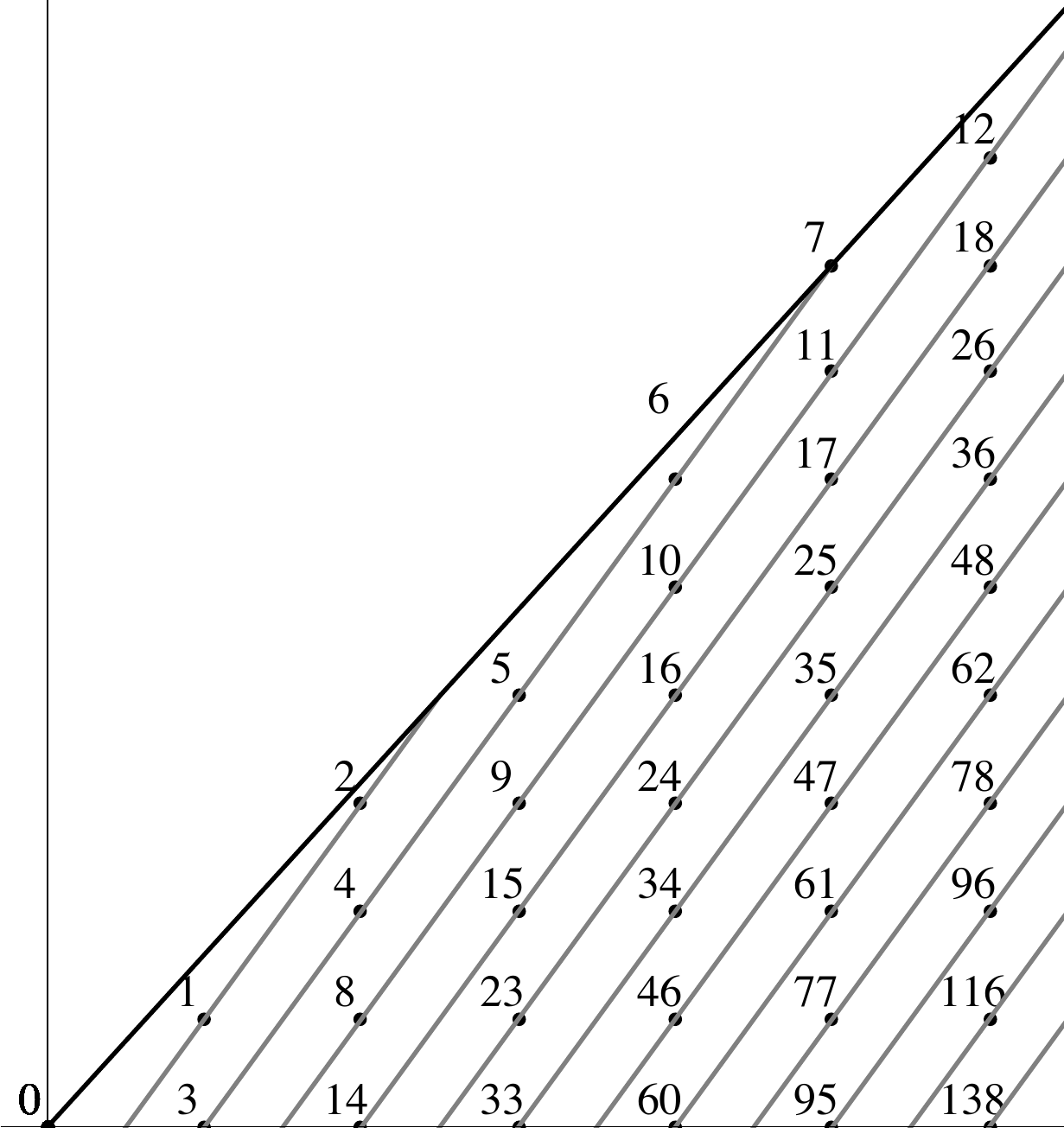}
  \caption{$p_+$}
  \label{fig:sub1}
\end{subfigure}%
\begin{subfigure}{.5\textwidth}
  \centering
  \includegraphics[width=.9\linewidth]{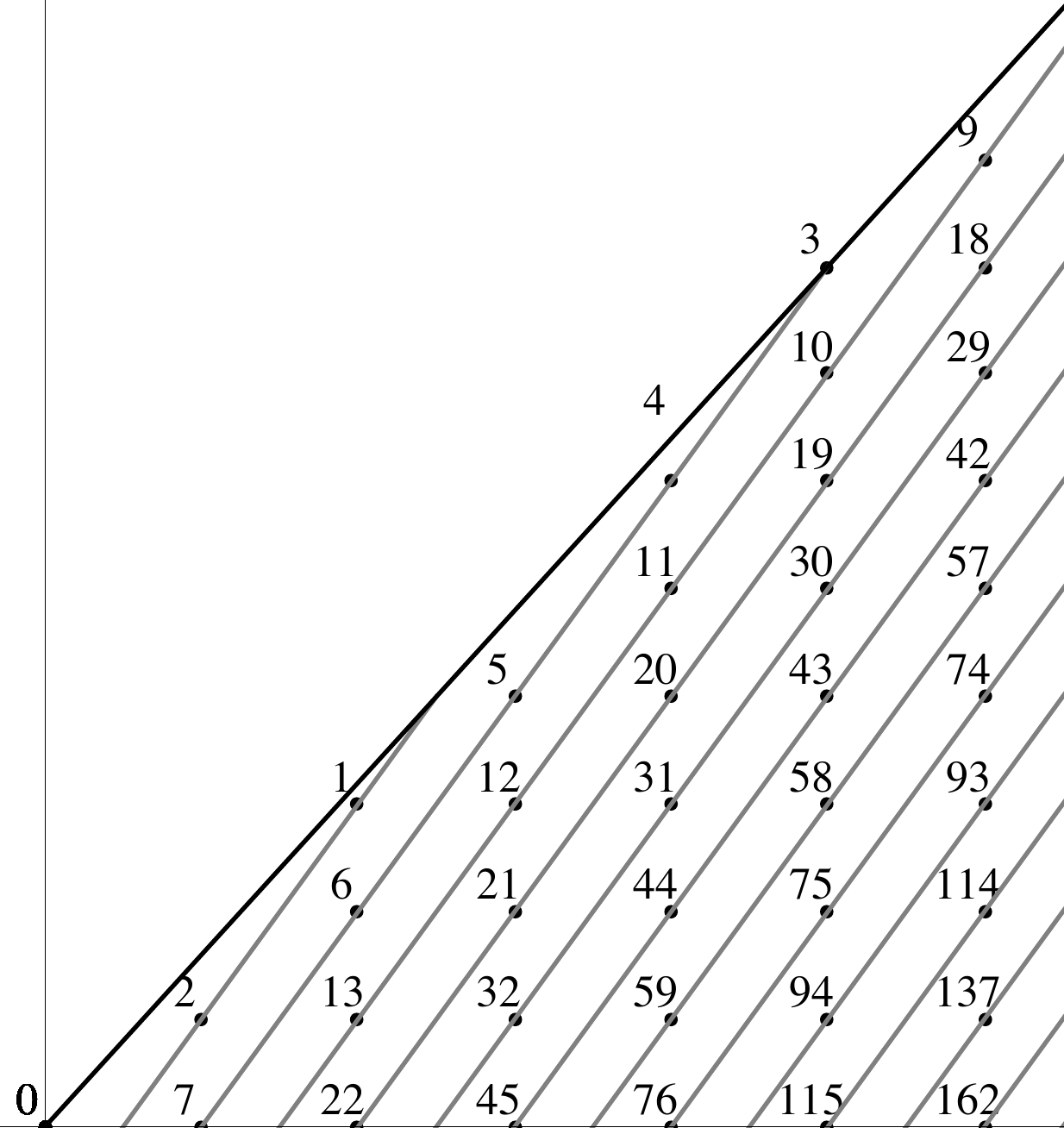}
  \caption{$p_-$}
  \label{fig:sub2}
\end{subfigure}
\caption{An ascending 1-stair packing polynomial $p_+$ and a descending 1-stair packing polynomial $p_-$, both on $S(\frac{8}{5})$.}
\label{4midn}
\label{fig:test}
\end{figure}

Lemma \ref{restriction1} shows that all quadratic packing polynomials on sectors $S(n/m)$ are $k$-stair polynomials for some $k$. Figure \ref{4midn} gives examples of two 1-stair packing polynomials. The next proposition shows that ascending and descending $k$ stair packing polynomials are equivalent.

\begin{proposition}
There is an ascending $k$-stair packing polynomial on $S(\frac{n}{m})$ if and only if there is a descending $k$-stair packing polynomial on $S(\frac{n}{n+2-m})$.
\begin{proof}
Let $m' = n+2-m$, and
$$
T_{n/m} = \left(
\begin{array}{c c}
m' & \frac{1-m'm}{n} \\
n & -m
\end{array}
\right).
$$
It is straightforward to show that $T_{n/m}$ is a bijection from $I(\frac{n}{m})$ to $I(\frac{n}{m'})$, so that
%%%% PROOF THAT IT IS A BIJECTION
%Since $\det(T) = -1$ and $T$ is an integer matrix, $T$ gives a bijection from $\ZZ^2$ to $\ZZ^2$. Let $(a,b) \in I(\frac{n}{m})$. Then
%$$
%T_{n/m}(a,b) = \left(m'a+\frac{1-m'm}{n}b , na-mb\right).
%$$
% Since $\frac{b}{a} \leq \frac{n}{m}$ and $a,b \geq 0$, we find that $0 \leq na-mb$. We will have that $0 \leq m'a+\frac{1-m'm}{n}b$ if and only if $0 \leq m'(na-mb)+b$, and since $0 \leq na-mb$, this holds. Then we must show that 
% $$
% \frac{na-mb}{m'a+\frac{1-m'm}{n}b} \leq \frac{n}{m'}.
% $$
%When $(a,b) \not =(0, 0)$, this inequality holds if and only if $0 \leq b$, which is true. Therefore, $T(a,b) \in I(\frac{n}{m'})$ and $T_{n/m} : I({\frac{n}{m}}) \rightarrow I(\frac{n}{m'})$ is injective.
%
%This implies that $T_{n/m'}: I({\frac{n}{m'}}) \rightarrow I(\frac{n}{m})$ is also injective, and since $T_{n/m}^{-1} = T_{n/m'}$, $T_{n/m} : I({\frac{n}{m}}) \rightarrow I(\frac{n}{m'})$ is bijective.\\
if $p$ is a quadratic packing polynomial on $S(\frac{n}{m})$, then $p \circ T_{n/m'}$ is a quadratic packing polynomial on $S(\frac{n}{m'})$. A simple calculation shows that if $p$ is an ascending (descending) $k$-stair polynomial, then $p \circ T_{n/m'}$ is a descending (ascending) $k$-stair packing polynomial on $S(\frac{n}{m'})$. 

%%%%% PROOF THAT IT IS KSTAIR
% Suppose $p$ is an ascending $k$-stair polynomial. Then for any $(x,y) \in I(\frac{n}{m})$,
%$$
%p\left(x+ \frac{m-1}{l},y+\frac{n}{l}\right) - p(x,y) = k.
%$$
%Also,
%\begin{align*}
%p\circ T_{n/m'}\left(x+ \frac{m-1}{l},y+\frac{n}{l}\right) - T_{n/m'}(x,y) &= p\left(T_{n/m'}(x+ \frac{m-1}{l},y+\frac{n}{l})\right) - p(T_{n/m'}(x,y)) \\
%& = p\left((x',y')+T_{n/m'}\left(\frac{m-1}{l},\frac{n}{l}\right)\right) - p(x',y')\\
%&= p\left(x'-\frac{m-1}{l}, y'-\frac{n}{l}\right) - p(x',y')\\
%&= -k.
%\end{align*}
%Therefore, $p \circ T_{n/m'}$ is a descending $k$-stair packing polynomial on $I(\frac{n}{m'})$. Similarly, if  $p$ is a descending $k$-stair polynomial on $I(\frac{n}{m})$, then
%$p \circ T_{n/m'}$ is an ascending $k$-stair packing polynomial on $I(\frac{n}{m'})$.
\end{proof}
\label{adthm}
\end{proposition}

%%%%%%%%%%%%%%
\subsection{Properties of $k$-stair polynomials.}
\begin{figure}
\begin{center}
\includegraphics[height = 2.6 in]{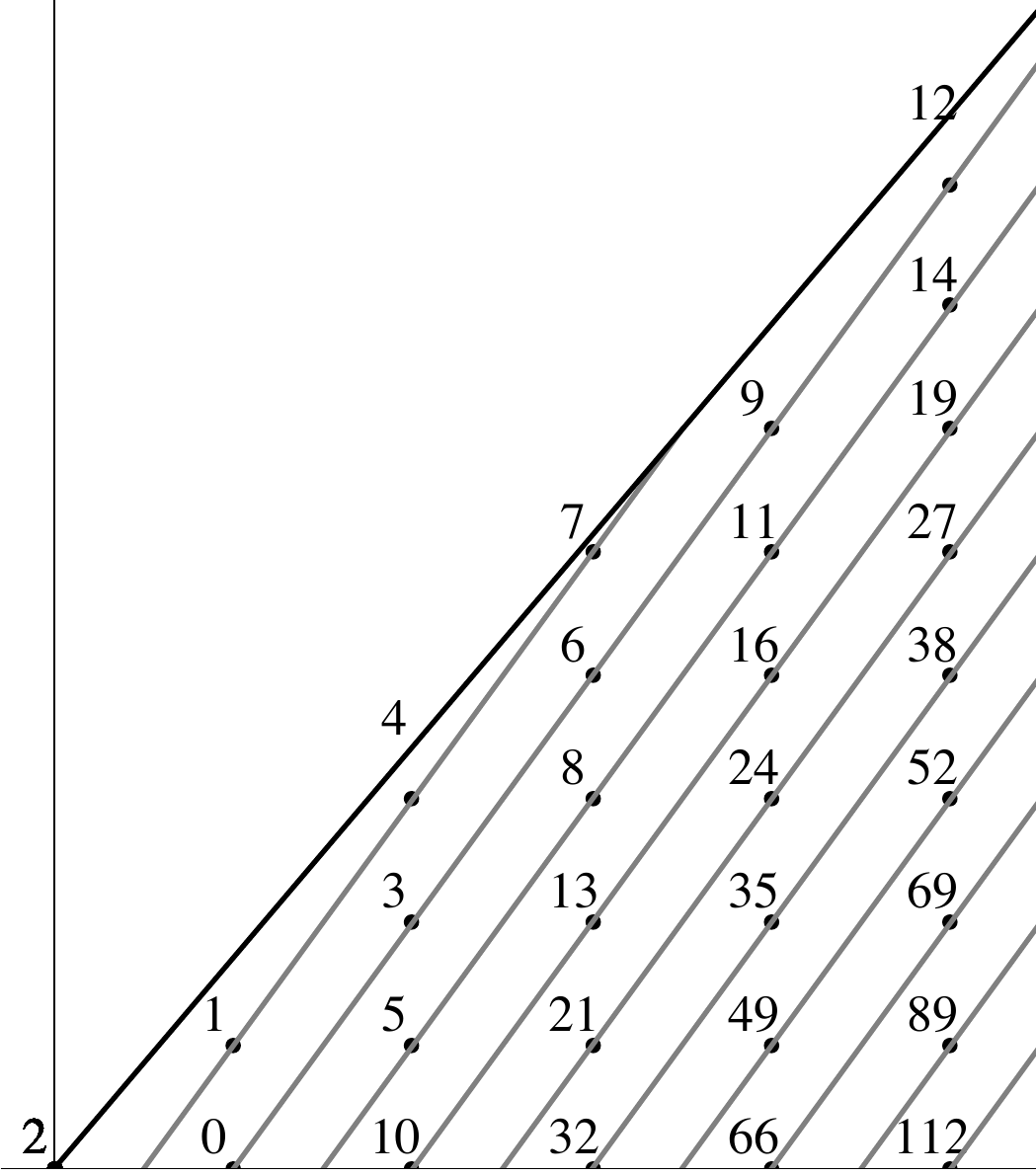}
\caption{An ascending 3-stair polynomial on $S(12/7)$.}
\label{tears}
\end{center}
\end{figure}

Let $p$ be a $k$-stair polynomial. Then the following immediate observations can be made.
If $a,b$ lie on the same staircase, then $p(a) \equiv p(b) \mod k$.
%\item $p(x,y) \equiv p(x+k\frac{l}{n},y) \mod k$.
%\item Let $a_1, \ldots, a_k$ be stairs from the $j^{th}$, $j+1^{st}, \ldots, j+k-1^{st}$ staircases for some $j$. Then $\{p(a_i) \mod k\ |\ i \in \{1, \ldots, k\}\} = \{1 , \ldots , k\}.$
Moreover, the numbers $0,1,\ldots, k-1$ must all occur on the first (last) stairs for an ascending (descending) $k$-stair packing polynomial, because otherwise the first (last) stairs will take on negative values.
%\item $p(0,0)+k = p(x_k,y_k)$ if $(x_k,y_k)$ is the first stair on the $k$th staircase.
%\item $f<k$.
Figure \ref{tears} gives an example of a 3-stair packing polynomial. The following lemma provides more information about the behavior of $k$-stair polynomials.

\begin{lemma} There exists some $j_0$ such that whenever $j\geq j_0$, if $(a,b) \in S_j$, and $(a',b') \in S_{j+k}$, then 
$$
p(a,b) \equiv p(a',b') \mod k.
$$
\begin{proof} The function $p(x,0)$ is increasing for all $x>x_0$ for some $x_0$. Let $j_0 = x_0 \frac{n}{l} $, and $j\geq j_0$. Suppose that for any $(a,b) \in S_j$, we have that $p(a,b) \equiv c \mod k$. Let $j'>j$ be the smallest integer such that for any $(a',b') \in S_{j'}$, we have that $p(a',b') \equiv c \mod k$. 

Suppose $j'-j > k$. Then by the pigeonhole principle, there exist $i,i'$ such that for any $(a,b) \in S_i$ and any $(a',b') \in S_{i'}$, $p(a,b) \equiv p(a',b') \mod k$. Let $s(l) $ be the number of stairs on the $l^{th}$ staircase, and let $\bar{p}(l)$ be the value of $p$ on the first stair on the $l$th staircase. Note that $s(i) = \lfloor(m-1)i \frac{l}{n}\rfloor \geq \lfloor (m-1)(j \frac{l}{n} + \frac{l}{n}) \rfloor \geq \lfloor(m-1)j\frac{l}{n}\rfloor +1$.

Then,
\begin{align*}
p\left(j' l/n,0\right) &\leq \bar{p}(j') \\
& = \bar{p}(j)+k \cdot s(j) \\
& < \bar{p}(i)+k \cdot s(j) \\
& \leq \bar{p}(i)+k \cdot (s(i)-1)\\
&= \bar{p}(i') - k\\
&<p\left(i' l/n,0\right),
\end{align*}
which is a contradiction because $p(x,0)$ is strictly increasing for $x > x_0$.
Therefore, $j'-j \leq k$. If $j'-j<k$, then there exist $i,i'$ such that  $i'-i > k$ and for any $(a,b) \in S_i$ and any $(a',b') \in S_{i'}$, we have $p(a,b) \equiv p(a',b') \mod k$, but we previously showed that this can not happen.
 Therefore, $j'-j=k$.
\label{lemmamodk}
\end{proof}
\end{lemma}

\begin{lemma}
Let $p$ be an ascending $k$-stair packing polynomial.
\begin{enumerate}
\item  If $j$ is a (large enough) integer, and $(a,b)$ is the first stair on the $(j\frac{n}{l}+k)^{th}$ staircase, then
$$
p(a,b)-k= p(mj,nj).
$$
\item If $(a,b)$ is the first stair on the $j^{th}$ staircase (for large enough $j$), and $(c,d)$ is the intersection of the line $y=\frac{n}{m}x$ and the $j-k^{th}$ staircase, then
$$
p(a,b)-k \leq p(c,d).
$$
\end{enumerate}

\begin{proof} \hfill

\begin{enumerate}
\item This is an immediate consequence of the bijectivity of $p$ and Lemma \ref{lemmamodk}.
\item If $p(a,b)-k > p(c,d)$ then the value $p(a,b)-k$ does not occur on any $i^{th}$ staircase where $i \equiv j \mod k$. However, if $p(x,y) \equiv p(a,b) \mod k$, and $p(x,y)$ is on the $i^{th}$ staircase, then $i \equiv j \mod k$, by Lemma \ref{lemmamodk}. Therefore, the value $p(a,b)-k$  is missing from the range of $p$, so $p$ is not surjective.
\end{enumerate}
\end{proof}
\label{restriction2}
\label{ineqlemma}
\end{lemma}

%%%%%%%%%
%%%%%%%%% n/l | k lemma
%%%%%%%%%
%\begin{lemma}[MB]
%If $\frac{n}{l} | k $ then there are no $k$-stair packing polynomials on $S(n/m)$.
%\label{nodivide}
%\begin{proof}
%Let $p$ be an ascending $k$-stair polynomial on $S(n/m)$, where $\frac{n}{l} | k$. From Lemma \ref{restriction1} and Lemma \ref{restriction2}, we find
%$$
%k = p\left(x+\frac{(m-1)}{l}, y+\frac{n}{l}\right)- p(x,y)
%$$
%and
%$$
%k= p\left(j+k\frac{l}{n},0\right) - p(m j, n j).
%$$
%From these restrictions, we find
%$$
%e= \frac{2n-2mn+kl^2+kml^2}{2nl}
%$$
%$$
%d=\frac{-n}{2}+\frac{n}{l}-\frac{kl}{2}.
%$$
%By Lemma \ref{ineqlemma}, we also require
%$$
%p\left(a-\frac{(m-1)}{l},b-\frac{n}{l}\right) \leq p\left(m\left(b\frac{(1-m)}{n} +a-k\frac{l}{n}\right),n\left(b\frac{(1-m)}{n} +a-k\frac{l}{n}\right)\right)
%$$
%will be satisfied when $(a,b)$ is a first stair. Simplifying, the inequality is satisfied if and only if
%$0 \leq k(\frac{l}{n} -1)$,
%which is never true because $l<n$. The descending case follows similarly.
%\end{proof}
%\end{lemma}

\begin{lemma}
Let $r$ be the multiplicative inverse of $\frac{m-1}{l} \mod \frac{n}{l}$. Let $j \equiv j' \mod \frac{n}{l}$ and $z = -j'r \mod \frac{n}{l}$. Then first stair on the $j^{th}$ staircase has coordinates
$$\left(\frac{m-1}{n}z+j\frac{l}{n},z \right).$$
\begin{proof}
Let $\phi : \ZZ_{\frac{n}{l}} \rightarrow  \ZZ_{\frac{n}{l}}$ send $j \in \ZZ_{\frac{n}{l}}$ to the $y$ coordinate of the first stair of the $j^{th}$ staircase. We only define $\phi$ on the first stairs because if $j \equiv i \mod \frac{n}{l}$, then the $y$ coordinate of the first stair of the $j^{th}$ staircase will be the same as the $y$ coordinate of the first stair of the $i^{th}$ staircase. Note that $\phi(j) \in \ZZ_{n/l}$ because if $(x,y)$ is the first stair on the $j^{th}$ staircase and $\phi(j) > \frac{n}{l}$, then $(x - \frac{m-1}{l}, y - \frac{n}{l})$ is a stair on the $j^{th}$ staircase with smaller $x$-coordinate.

 Then observe that $\phi^{-1}(-i) = i \frac{m-1}{l}$, so that $\phi(j)=-j(\frac{m-1}{l})^{-1}$. The $x$ coordinate comes from solving
$$
y = \frac{n}{m-1}\left(x - j \frac{l}{n}\right).
$$ 
\end{proof}
\label{stairlemma}
\end{lemma}

%%%%%%%%%%%%%%%%
\section{Necessary and Sufficient Conditions for $k$-Stair Packing Polynomials on $S(\frac{n}{m})$}
\label{nsconditions}

\begin{theorem}
Let $p(x,y) = \frac{n}{2} (x - \frac{m-1}{n} y)^2+dx+ey+f$ be a packing polynomial on $S(n/m)$, where $l = (n,m-1)$. Then either
\begin{enumerate}
\item $p$ is an ascending $k$-stair polynomial where $k\equiv \frac{m-1}{l}\mod\frac{n}{l}$, and
$$
p(x,y) = \frac{n}{2}\left(x-\frac{m-1}{n}y\right)^2+\left(1-\frac{kl}{2}\right)x+
\frac{2(1-m)+kl(m+1)}{2n}y+f,
$$
or
\item$p$ is a descending $k$-stair polynomial where $k \equiv -\frac{m-1}{l}\mod\frac{n}{l}$, and
$$
p(x,y) = \frac{n}{2}\left(x-\frac{m-1}{n}y\right)^2+
\left(1+\frac{kl}{2}\right)x+
\frac{2(1-m)-kl(m+1)}{2n}y
+f.
$$
\end{enumerate}
\begin{proof}
Let $p$ be an ascending $k$-stair packing polynomial on $S(n/m)$. Let $k = q\frac{n}{l}+k'$, where $k'<\frac{n}{l}$ and $q$ is an integer; let $r$ be the multiplicative inverse of $\frac{m-1}{l} \mod \frac{n}{l}$, and let $z = -k'r \mod \frac{n}{l}$. By Lemma \ref{restriction1} and Lemma \ref{restriction2}, for large enough $j$ and $x$, we have
$$
k = p\left(x+\frac{(m-1)}{l}, y+\frac{n}{l}\right)- p(x,y)
$$
and
$$
k= p\left(\frac{m-1}{n}z+j+k\frac{l}{n},z \right) - p(m j, n j),
$$
so that
$$
e= \frac{-2(m-1)n+((m+1)nq+2(m-1)z)l+k'(m+1)l^2}{2nl},
$$
and
$$
d=\frac{1}{2}\left(-nq-2z+\frac{2n}{l}-k'l\right).
$$
By Lemma \ref{ineqlemma}, we also have when $(a,b)$ is a first stair (for large enough $a$),
\begin{align*}
p\left(a-\frac{(m-1)}{l},b-\frac{n}{l}\right) \leq p\Big(&m\left(b\frac{(1-m)}{n} +a-(q\frac{n}{l}+k')\frac{l}{n}\right),\\
& n\left(b\frac{(1-m)}{n}+a-(q\frac{n}{l}+k')\frac{l}{n}\right)\Big).
\end{align*}
Plugging in these points using the $e$ and $d$ given above, we find that this inequality is satisfied if and only if
$$
0 \leq -\frac{(b-z)(nq+k'l)}{n}.
$$
In particular, there are first stairs $(a,b)$ where $b = \frac{n}{l}-1$, and because $\frac{nq+k'l}{n}> 0$, we therefore must have that $0 \leq z-b$. Since $z < \frac{n}{l}$, this implies that $z = \frac{n}{l}-1$. So, $-k'r \equiv -1$, so $k' \equiv r^{-1}$, which implies that $k' = \frac{m-1}{l} \mod \frac{n}{l}$. The coefficients of $p(x,y)$ follow from simplifying $d$ and $e$ with this requirement.\\

The case where $p$ is a descending polynomial follows from Proposition \ref{adthm}.
\end{proof}
\label{bigthm}
\end{theorem}

Since the coefficients $e,f$ must be set to satisfy the inequalities from Lemma \ref{ineqlemma} for large enough $x,y$, they will also satisfy the inequalities for all other $x,y$. Therefore, we find that $p(x,y)$ automatically satisfies the inequalities from Lemma \ref{ineqlemma} for all $x,y$.

\begin{theorem}
Let $a_1, \ldots, a_{k}$ be the first stairs on the first $k$ staircases on $S(n/m)$. Then $p$ is an ascending packing polynomial if and only if $p$ is a $k$-stair polynomial of the necessary form given in Theorem \ref{bigthm}, and
$$
\{p(a_1), \ldots, p(a_{k})\} = \{0,1, \ldots, k-1\}.
$$
\begin{proof}
Suppose $p$ is an ascending $k$-stair polynomial with the necessary form, and $\{p(a_1), \ldots, p(a_{k})\} = \{0,1, \ldots, k-1\}$. For any $i \in \{1, \ldots, k\}$, let 
$$R_i =\cup\{S_c\ |\ c \equiv i \mod k\}.$$
If $p|_{R_i}$ is a bijection from  $R_i\cap \NN^2$ to $p(a_i)+k \NN$ for any $i$, then $p$ is a packing polynomial on $S(n/m)$.  

 Since $p$ satisfies the inequality from Lemma $\ref{ineqlemma}$, $p|_{R_i}$ is surjective to $p(a_i)+k \NN$ (since no values congruent to $i \mod k$ will be skipped). Then $p|_{R_i}$ will be injective if whenever $(a,b)$ is the first stair on the $j^{th}$ staircase and $(c,d)$ is the last stair on the $(j-k)^{th}$ staircase, we have
 $$
 p(a,b)>p(c,d).
 $$
We also have
 $$
0< \frac{m-1+nq}{n} = p\left(j\frac{l}{n},0\right)-p\left(m\frac{l}{n}(j-k),l(j-k)\right).
 $$
Since 
 $$
 p(a,b) \geq p\left(j\frac{l}{n},0\right)> p\left(m\frac{l}{n}(j-k),l(j-k)\right) \geq p(c,d),
 $$
 we have $p(a,b)>p(c,d)$, so $p|_{R_i}$ is injective and so $p$ is a packing polynomial.\\

Conversely, suppose $p$ is a packing polynomial. Let $i \in \{0,1,\ldots,k-1\}$. If $p(a_i) \equiv j \mod k$ where $0 \leq j < k$, but $p(a_i) \not = j$, then by the above, for any $a \in R_i$, we have $p(a) \geq p(a_i)$. So, there is no $(x,y)$ such that $p(x,y) = j$. Therefore, $p(a_i) = j$ for some $j \in \{0,1, \ldots, k-1\}$. On the other hand, if there is some $j \in \{0,1,\ldots, k-1\}$ such that $p(a_i) \not = j$ for any $i \in \{1,\ldots,k\}$, then $p$ will never achieve values congruent to $j \mod k$. Therefore,  $\{p(a_1), \ldots , p(a_{k})\} = \{0,1, \ldots, k-1\}$.\\

\end{proof}
\label{sufficient}
\end{theorem}

%%%%%%%%%%%%%%
\section{Packing Polynomials on $S(\frac{n}{m})$}
\label{knot4}

Now we are prepared to determine, up to isomorphism, the $k$-stair packing polynomials for each $k$. In particular, we will prove that there are no $k$-stair polynomials when $k \geq 4$.
We first provide two additional results.

\begin{proposition} If there is an ascending $k$-stair packing polynomial on $S(\frac{n}{m})$ for $\frac{n}{m}> 1$ and $m\not=1$, then $k \mid l$.

\begin{proof}
Suppose $l \not \equiv 0 \mod k$. Let $(a_i, b_i)$ be the first stair on the $i^{th}$ staircase, where $i < k$. Let $\xi_i = - i (\frac{m-1}{l})^{-1} \mod \frac{n}{l}$. Then by Lemma \ref{stairlemma} and Theorem \ref{bigthm} we have
$$
p(a_i,b_i) - f = \frac{l}{2n}(i (i-k) l + 2(i + k \xi_i)). 
$$
By Theorem \ref{sufficient}, for any $i , j < k$, we need that $p(a_i,b_i) \not \equiv p(a_j,b_j) \mod k$. Since $\frac{n}{l}$ and $k = q \frac{n}{l} + \frac{m-1}{l}$ (for some integer $q$) are relatively prime, $n/l$ is not a zero divisor in $\ZZ_k$. Therefore, $p(a_i,b_i) \not \equiv p(a_j,b_j) \mod k$ if and only if $\frac{n}{l}p(a_i,b_i) \not \equiv \frac{n}{l}p(a_j,b_j) \mod k$. Then,
\begin{align*}
\frac{n}{l}p(a_i,b_i) &= \frac{l }{2} i^2 + i\left(1 - \frac{k l}{2}\right)\\
&= \frac{l}{2}\left(i + \frac{1}{l} - \frac{k}{2}\right)^2-\frac{l}{2}\left(\frac{k}{2}-\frac{1}{l}\right)^2.
\end{align*}
So, if $j = -i -\frac{2}{l} +k$, then $p(a_i,b_i) \equiv p(a_j,b_j)\mod k$. If $j = i$, then $i = \frac{k}{l} - \frac{1}{l}$, which will only happen for one $i$, so for some $i,j$, we have that $p(a_i,b_i) \equiv p(a_j,b_j)\mod k$.
\end{proof}
\label{lmodkthm}
\end{proposition}

%%%%% Computation
%Proposition \ref{lmodkthm} implies that $k \mid l$, so that for a given sector $S(n/m)$, there are only finitely many possible $k$ values. Using this observation, a computer search reveals that for $n,m \leq 50,000$, no $k$-stair packing polynomials exist on these sectors for $k \geq 4$. This leads us to make the following conjecture.
%
%\begin{conjecture} There are no $k$ stair packing polynomials for $k \geq 4$.
%\label{4pstair}
%\end{conjecture}
%
%As further evidence for Conjecture \ref{4pstair}, we have the following result.

\begin{theorem}
If $p$ is a $k$-stair packing polynomial on $S(\frac{n}{m})$ where $\frac{n}{m}>1$ and $m\not=1$, then either $k = \frac{m-1}{l}$ or $\frac{n}{m} = \frac{12}{7}$.
\begin{proof}
By Theorem \ref{bigthm}, we have $k = q\frac{n}{l} + \frac{m-1}{l}$ for some $q \in \NN$. Suppose that $q \not = 0$. Then the point $(1,0)$ is the first stair on the $\frac{n}{l}^{th}$ staircase, and substituting $l = q\frac{n}{k} + \frac{m-1}{k}$, we find that
$$
p(1,0) -f= \frac{1}{2}(3-m+n-nq).
$$
By Theorem \ref{sufficient}, we have $|p(1,0)-f|\leq k-1$. In particular, the inequality $p(1,0) - f \geq -(k-1)$ implies that
$$
2k \geq m-1+n(q-1).
$$
If $q>1$, then 
\begin{align*}
2k &\geq m-1+n(q-1)\\
&\geq m-1+n\\
&>2(m-1),
\end{align*}
which implies that $k > m-1$. This is impossible; by Proposition \ref{lmodkthm}, we have $k \mid l$, and $l \mid m-1$, so $k \leq m-1$ (when $m \not = 0$).\\

Suppose that $q=1$. Then $2k \geq m-1 \geq k$, and since $k \mid m-1$, we must have $2k = m-1$, or $k = m-1 = l$. 
\begin{enumerate}
%%%%%%%%%%
\item Suppose that $2k=m-1$. Then $k \mid l$ and $l \mid 2k$, so either $l = 2k = m-1$ or $l = k$.

\begin{enumerate}
%%%%%%%%%% 	1A
\item Suppose that $l = 2k = m-1$. Then the equality $k = q \frac{n}{l} + \frac{m-1}{l}$ implies that 
%$k = \frac{1 + \sqrt{1+2n}}{2}$, and $m = 2+ \sqrt{1+2n}$.
$k = \frac{n}{m-1}+1$, so 
$$
\frac{(m-1)^2}{n} = \frac{m-1}{k-1} = \frac{2k}{k-1}.
$$
Moreover, this quantity is an integer. We conclude that either $k =2$ or $(k - 1) \mid 2$, so that $k = 2$ or $3$.

If $k =2$, then $m=5$ and $n=4$, contradicting the assumption that $n>m$.

If $k =3$, then $m = 7$ and $n = 12$, and we obtain a 3-stair packing polynomial on $S(\frac{12}{7})$.

%%%%%%%%%%		1B
\item Suppose that $2l = 2k =m-1$. Then the equality $k = q \frac{n}{l} + \frac{m-1}{l}$ implies that 
%$k = 1+\sqrt{1+n}$, and $m = 2+ \sqrt{1+2n}$.
$k = \frac{2n}{m-1}+2$, so
$$
\frac{(m-1)^2}{n} = \frac{2(m-1)}{k-2} = \frac{4k}{k-2}.
$$
Moreover, this quantity is an integer. For any integer $k$, $\gcd(k,k-2) = 1$ or 2. 

If $\gcd(k,k-2) = 1$, then $(k-2) \mid 4$, which forces $k=3$. If $k = 3$, then $n = 3$ and $m = 7$, contradicting the assumption that $n>m$.

If $\gcd(k,k-2) = 2$, then $2 \mid k$, so $\frac{n}{m-1}$ is an integer, forcing $l = m-1 = 2l$, so this case cannot occur.

%$\frac{k-2}{2} \mid 4$, so $k = 10, 6,$ or 4. Therefore, $k$ could be 10, 6, 4, or 3.
%
%If $k=10$, then $n=10$ and $m = 21$, and so $l = 20 \not = k$, which is a contradiction.
%
%If $k = 6$, then $n=24$, and $m = 13$, and so $l = 12 \not = k$, which is a contradiction.
%
%If $k = 4$, then $n=8$ and $m = 9$, which is a contradiction because $n>m$.

\end{enumerate}
%%%%%%%%%%		2
\item Suppose that $k = m-1 = l$. Then the equality $k = q \frac{n}{l} + \frac{m-1}{l}$ implies that
% $k =\frac{1+\sqrt{1+4n}}{2}$, and $m = \frac{3+\sqrt{1+4n}}{2}$.
 $k = \frac{n}{m-1}+1$, so
$$
\frac{(m-1)^2}{n} = \frac{m-1}{k-1} = \frac{k}{k-1}.
$$
Moreover, this quantity is an integer, so that $k =2$. It follows that $m=3$ and $n = 2$, contradicting the assumption that $n>m$.\\
\end{enumerate}

Therefore, if $k \not = \frac{m-1}{n}$, then $\frac{n}{m} = \frac{12}{7}$.

\end{proof}
\label{qzero}
\end{theorem}

%%%%%%%%%%

\begin{theorem} Let $\frac{n}{m} \in \QQ$, $(n,m) = 1$, $m \not = 1$, and $\frac{n}{m} > 1$. The following results give the $k$-stair packing polynomials on sectors $S(\frac{n}{m})$ for $k \in \{1,2,3,4\}$.
\begin{enumerate}
\item  There is an ascending 1-stair packing polynomial $p$ on $S(n/m)$ if and only if  $n \mid (m-1)^2$ and $m-1 \mid n$.
\item There is an ascending 2-stair packing polynomial $p$ on $S(n/m)$ if and only if $m \equiv 9 \mod 16$ and $n = \frac{1}{16}(m-1)^2$.
\item
There is an ascending $3$-stair packing polynomial $p$ on $S(n/m)$ if and only if  $m \equiv 10 \mod 27$ or $m \equiv 19 \mod 27$ and $n = \frac{1}{27}(m-1)^2$, or $\frac{n}{m} = \frac{12}{7}$.
\item There are no $4$-stair packing polynomials.
\end{enumerate}

\begin{proof}
\begin{enumerate}
%%%%%%%%% 			k=1
\item By Theorem \ref{bigthm}, we have $k=1$ if and only if $\frac{m-1}{l} = 1$, so $m-1 \mid n$.
 By setting $f=0$, the sufficient condition from Theorem \ref{sufficient} is satisfied since $p(0,0) = f=0$.
%%%%%%%%% 			k=2
\item Suppose $p$ is an ascending 2-stair polynomial on $S(n/m)$.
By Theorem \ref{qzero}, $k = \frac{m-1}{l} = 2$, $2 \nmid \frac{n}{l}$. Note that $2(\frac{n}{l} - \frac{n/l-1}{2}) \equiv 1 \mod \frac{n}{l}$, so the first stair on the first staircase (by Lemma \ref{stairlemma}) is
$$
\left(1, \frac{n/l-1}{2}\right).
$$
Then by Theorem \ref{sufficient}, $p$ is a packing polynomial if and only if 
$$
\left\{p(0,0), p\left(1,\frac{n/l-1}{2}\right)\right\} = \{0,1\}.
$$
Since $p(0,0) = f$, we may find an $f$ which satisfies this as long as $|p(0,0)- p(1,\frac{n/l-1}{2})| = 1$. So, using the necessary form of $p$ given in Theorem \ref{bigthm} we have,
\begin{align*}
\pm 1 &= p\left(1,\frac{n/l-1}{2}\right)-p(0,0)\\
&=  \frac{(m-1)^2 - 8n}{8n}.
\end{align*}
Since $m\not=1$, we find that $p$ is a packing polynomial if and only if $n= \frac{(m-1)^2}{16}$. Then because $l = \frac{m-1}{2}$, we have that $8 \mid m-1$ but $16 \nmid m-1$.

In the case where $m=1$, Stanton in \cite{stanton} found two 2-stair packing polynomials on $S(4)$. However, by the above we note that $S(4/9)$ has a 2-stair packing polynomial, and Stanton's polynomials are both equivalent to the ascending 2-stair packing polynomial on $S(4/9)$.

%%%%%%%%% 			k=3
\item By Theorem \ref{qzero}, either $n/m = 12/7$, or $k = 3 = \frac{m-1}{l}$. The case where $3 = \frac{m-1}{l}$ follows by a method similar to that used to prove (2). We note that the ascending 3-stair packing polynomial on $S(12/19)$ is equivalent to the ascending 3-stair packing polynomial on $S(12/7)$.

In the case where $m=1$, Stanton in \cite{stanton} found two 3-stair packing polynomials on $S(3)$. In a fashion similar to $(2)$, these are both equivalent to the ascending 3-stair packing polynomial on $S(3/10)$.\\

%%%% THE OLD PROOF OF THIS, BEFORE THE QZERO THEOREM
%Since $k = 3 \equiv \frac{m-1}{l} \mod \frac{n}{l}$, either $\frac{m-1}{l} = 3$ or $\frac{n}{l}=2$ and $\frac{m-1}{l} = 1$.
%\begin{enumerate}
%\item If $\frac{m-1}{l} = 3$, then (1) follows by a method similar to that used to prove Theorem \ref{k2}.
%\item If $\frac{n}{l}=2$ and $\frac{m-1}{l} = 1$, then by Lemma \ref{stairlemma}, the first stairs of the first 3 staircases are $(0,0), (1,1),$ and $(1,0)$. Since $p(0,0)=f$, we need $\{p(1,1)-f,p(1,0)-f\} \subset \{-2,-1,1,2\}$. Using the necessary equation for $p$ from Theorem \ref{bigthm},
%$$
%p(1,1)-f = \frac{3-m}{2} \ \ \ \ \ \ \ \ \ p(1,0)-f = \frac{5-m}{2},
%$$
%so $p(1,1) = p(1,0)-1$. Since $m \not = \pm 1$ and $p(1,1)-f \not = 1,2$. If $p(1,1)-f = -1$, then $p(1,0)-f = 0$, which can not happen. Then $p(1,1)-f = -2$, $p(1,0)-f = -1$, and $m=7$ and $n=12$. Setting $f=2$, the necessary $p$ is a 3-stair ascending packing polynomial on $S(12/7)$ by Theorem \ref{sufficient}.
%\end{enumerate}
%%%%%%%%% 			k=4
\item Stanton \cite{stanton} proved that there are no 4-stair packing polynomials on $S(n)$.
By Theorem \ref{qzero}, we have that if $m\not=1$, then $4 =\frac{m-1}{l}$.
Since $\frac{n}{l}$ is relatively prime to $\frac{m-1}{l} = 4$, either $\frac{n}{l}\equiv 1 \mod k$ or $\frac{n}{l}\equiv 3 \mod k$.
\begin{enumerate}
\item Suppose $\frac{n}{l}\equiv 1 \mod k$.
Then if $(a_i,b_i)$ are the first stairs on the first $4$ staircases by Lemma \ref{stairlemma} we have
\begin{align*}
(a_0,b_0) &= (0,0) \\
(a_1,b_1) &= \left(1, \frac{n/l-1}{4}\right)\\
(a_2,b_2) &= \left(2, \frac{n/l-1}{2}\right)\\
(a_3,b_3) &= \left(3, 3\frac{n/l-1}{4}\right),
\end{align*}
and so 
\begin{align*}
p(a_0,b_0)-f &=0 \\
p(a_1,b_1)-f &= \frac{-3 (m-1)^2+32n}{32 n}\\
p(a_2,b_2) -f&= \frac{(m-1)^2 - 16n}{8n}\\
p(a_3,b_3) -f&= \frac{-3((m-1)^2- 32 n)}{32n}.
\end{align*}
By Theorem \ref{sufficient}, $\{p(a_i,b_i)-f \ |\ i \in \{0,1,2,3\}\} \subset \{-3,-2,-1,-,1,2,3\}$.

Let $x = p(a_1,b_1)-f $. Then $32 n (1-x) = 3(m-1)^2$, so $x<0$. Let $x' = p(a_2,b_2) -f$. Then $-8nx'+16n=(m-1)^2$, so
$$
8n(2-x') =\frac{ 32 n}{3}(1-x), 
$$
which implies that $3 \mid 1-x$ and $4 \mid 2-x'$, so $x=x' = -2$. We conclude that $p$ is not injective.
%%%%%%%%
\item The case where $\frac{n}{l}\equiv 3 \mod k$ follows similarly; we find that there are no packing polynomials in this case.

%%% PROOF OF $\frac{n}{l}\equiv 3 \mod k$ CASE
%%%
%%%
%Suppose $\frac{n}{l}\equiv 3 \mod k$.
%Then if $(a_i,b_i)$ are the first stairs on the first $4$ staircases by Lemma \ref{stairlemma} we have
%\begin{align*}
%(a_0,b_0) &= (0,0) \\
%(a_1,b_1) &= \left(3, \frac{3n/l-1}{4}\right)\\
%(a_2,b_2) &= \left(2, \frac{n/l-1}{2}\right)\\
%(a_3,b_3) &= \left(1, \frac{n/l-3}{4}\right),
%\end{align*}
%and so 
%\begin{align*}
%p(a_0,b_0)-f &=0 \\
%p(a_1,b_1)-f &= \frac{-3((m-1)^2-32n)}{32 n}\\
%p(a_2,b_2) -f&= \frac{(m-1)^2 - 16n}{8n}\\
%p(a_3,b_3) -f&= \frac{-3(m-1)^2+32n}{32n}.
%\end{align*}
%By Theorem \ref{sufficient}, $\{p(a_i,b_i)-f \ |\ i \in \{0,1,2,3\}\} \subset \{-3,-2,-1,-,1,2,3\}$.
%
%Let $x = p(a_1,b_1)-f $. Then $\frac{32 n (3-x)}{3} = (m-1)^2$. Let $x' = p(a_2,b_2) -f$. Then $8n(x'+2)=(m-1)^2$, so
%$$
%\frac{32 n (3-x)}{3}=8n(x'+2),
%$$
%which implies that $3 \mid 3-x$ and $4 \mid 2+x'$, so $x=0$, $x'=2$, and $p(a_3,b_3) -f = -2$, which contradicts Theorem \ref{sufficient}.

\end{enumerate}

%%%%%%%% THE EXTRA CASE FROM Q=1
%\item Suppose  $q=1$, $\frac{n}{l} = 3$, and $\frac{m-1}{l} = 1$.
%Then if $(a_i,b_i)$ are the first stairs on the first $4$ staircases by Lemma \ref{stairlemma},
%\begin{align*}
%(a_0,b_0) &= (0,0) \\
%(a_1,b_1) &=(1,2)\\
%(a_2,b_2) &=(1,1)\\
%(a_3,b_3) &= (1,0),
%\end{align*}
%and so 
%\begin{align*}
%p(a_0,b_0)-f &=0 \\
%p(a_1,b_1)-f &=\frac{7-m}{2}\\
%p(a_2,b_2) -f&= \frac{-2(m-4)}{3}\\
%p(a_3,b_3) -f&= \frac{3-m}{2}.
%\end{align*}
%So, $p(1,0)+2=p(1,2)$, which implies that $p(1,2)-f \not = 2,0,-1,-2$. If $p(1,2)-f = 3$, then $p(1,1)-f=2$, in which case $m=1$. If $p(1,2)-f = 1$, then $m=5$, so $p(1,1) = -2/3$, neither of which can happen.

\end{enumerate}
\label{1234thm}
\end{proof}
\end{theorem}
%%%%%%

\begin{theorem} There are no $k$-stair polynomials for $k \geq 4$.
\begin{proof} Stanton \cite{stanton} showed that in the case where $m=1$, there are no $k$-stair packing polynomials for $k \geq 4$, so assume $m \not = 1$.
Let $p$ be a $k$-stair packing polynomial where $k \geq 4$. By Theorem \ref{qzero}, this implies that $k = \frac{m-1}{l}$. Then by Theorem \ref{bigthm},
$$
p(x,y) = \frac{n}{2} \left(x - \frac{m-1}{n} y  \right)^2 + \frac{3-m}{2} x + \frac{(m-1)^2}{2n} y + f.
$$
Let $(a,b)$ be the first stair on $S_1$. By Theorem \ref{sufficient}, we know that 
$$|p(a,b) - f| \leq k-1.$$
 Also, $p(a,b) - p(\frac{l}{n}, 0) \leq k-1$ since $p$ is $k$-stair, so that $p(\frac{l}{n},0) - f \geq -2(k-1)$. From the above form of $p(x,y)$, we find that
$$
p\left(\frac{l}{n},0\right) = \frac{(m-1)^2}{2kn} \left(\frac{1-k}{k} + \frac{2}{m-1} \right).
$$
Then $p(\frac{l}{n},0) - f \geq -2(k-1)$ if and only if
$$
\frac{(m-1)^2}{kn} \left(\frac{1}{k} - \frac{2}{(m-1)(k-1)} \right) \leq 4.
$$
Now, we claim that $n \leq \frac{l^2}{k}$. Suppose that $n$ and $m-1$ are divisible by $k^j$, and no higher power of $k$ divides $m-1$. Since $m-1 = k l$, we have that $k^{j-1}$ is the highest power of $k$ that divides $l$. But, since $n$ and $m-1$ are both divisible by $k^j$, we also have $k^j \mid l$, which is a contradiction. Therefore, a higher power of $k$ divides $m-1$ than divides $n$. Also, $n \mid (m-1)^2$, so $n \mid l^2$. If $k^{i}$ is the highest power of $k$ that divides $l$ and $k^{2i} \mid n$, then by the above, $2i < i+1$, which implies that $i = 0$. This is a contradiction since $k \mid l$. Therefore, $n \mid \frac{l^2}{k}$.\\

Plugging in $n = \frac{l^2}{k}$, we find that
\begin{equation}
k - \frac{2k}{l(k-1)} \leq \frac{(m-1)^2}{kn} \left(\frac{1}{k} - \frac{2}{(m-1)(k-1)} \right) \leq 4.
\label{ineq}
\end{equation}
Therefore, if
$$
k > 4+ \frac{2k}{l(k-1)},
$$
then inequality (\ref{ineq}) does not hold, and so $p$ is not a packing polynomial. We find that $
k > 4+ \frac{2k}{l(k-1)}$ whenever $k \geq 5$. This result, along with Theorem \ref{1234thm}, implies that there are no $k$-stair packing polynomials when $k \geq 4$.
\end{proof}
\end{theorem}

Thus, there are only $k$-stair polynomials for $k \in \{1,2,3\}$. By Theorem \ref{1234thm} and extra details provided in the proof, we conclude that up to equivalence, the polynomials given in Theorem \ref{bigthm}.1 on the sectors $S(\frac{n}{m})$ given by Theorem \ref{1234thm} along with Nathanson's $f_n, g_n$ on $S(n)$ represent all quadratic packing polynomials on rational sectors.

\section{Future Directions}
Lew and Rosenberg \cite{LR1} proved that there are no packing polynomials of degree three or four on $\NN^2$. It is an open question whether there exist packing polynomials of degree greater than two on rational sectors. In addition, the conjecture of  Nathanson \cite{Nathonson} that there are no packing polynomials on $S(\alpha)$ for irrational $\alpha$ remains open.

\section{Acknowledgements}
I conducted this research while I was a participant of the University of Minnesota Duluth Research Experience for Undergraduates in Mathematics Program, supported by the National Science Foundation (grant number DMS 1358659) and the National Security Agency (NSA grant H98230-13-1-0273). I would like to acknowledge the valuable suggestions and encouragements of the program's visitors and advisors, especially Sam Elder and Sam Stewart. I would like to extend a special thanks to Joe Gallian for his enthusiasm and advice.

\nocite{*}
\bibliographystyle{plain}
\bibliography{packingpoly}

\end{document}